# Generalization of Trigonometric B-splines and Kernels of Interpolating Trigonometric Splines with Riemann Convergence Multipliers

**DENYSIUK Volodymyr Petrovych,**

**HRYSHKO Olena Mykolayivna**

**Abstract**

*This paper explores the generalization of the method for extracting Riemann trigonometric B-splines and Riemann kernels of trigonometric interpolation splines of arbitrary order on different grids of stitching and interpolation. It is demonstrated that for various combinations of stitching and interpolation grids, distinct trigonometric B-splines exist. The theoretical principles are illustrated through a numerical example. The obtained results can have various practical applications.*

## Introduction

Approximation, representing, any known or unknown function through a set of certain special functions, can be considered a central theme in analysis. The term 'special functions' refers to classes of algebraic and trigonometric polynomials and their modifications, with the classes of trigonometric polynomials including trigonometric series. Typically, such special functions are easily computable and possess interesting analytical properties [1].

One of the most successful modifications of algebraic polynomials is polynomial splines, which are stitched together from segments of these polynomials according to a certain scheme. The theory of polynomial splines appeared relatively recently and is well developed (see, e.g., [2], [3], [4], [5], etc.). The advantages of polynomial splines include their approximative properties [6]. The main disadvantage of polynomial splines is their piecewise structure, which significantly complicates their application in analytical transformations.

Normalized B-splines (simply B-splines) play an important role in the theory of polynomial splines, representing a basis in the space of polynomial splines [3], [4]. The representation of polynomial splines through B-splines is used in many problems requiring analytical transformations, particularly in solving differential and integral equations (see, e.g., [7], [8], etc.).

Later, it was discovered [9], [10], that there are modifications of trigonometric series that depend on several parameters and have the same properties as polynomial splines [11] [12]; moreover, the class of such modified series is quite broad and includes the class of polynomial periodic splines. This provided grounds to call the class of such series trigonometric splines.

Trigonometric splines and their generalizations were also considered in [13], [14]. In [15], a method for extracting trigonometric B-splines and kernels of interpolation trigonometric splines with Riemann convergence factors was proposed for one type of combination of stitching and interpolation grids. It is of great interest, in our opinion, to generalize this method to other types of combinations of these grids.

**The goal of the work.** Generalization of the method for constructing trigonometric B-splines and kernels of interpolation trigonometric splines with Riemann convergence factors, for various types of combinations of stitching and interpolation grids of trigonometric splines.

## Main part

Let a continuous periodic function $f(t)$ be defined on the interval $[0, 2\pi]$. Additionally, let a uniform grid $\Delta_N^{(I)} = \{t_i^{(I)}\}_{i=1}^{N}$, ($I = 0,1$) $t_i^{(0)} = \frac{2\pi}{N}(i-1)$, $t_i^{(1)} = \frac{\pi}{N}(2i-1)$, $N = 2n+1$, $n = 1, 2, ...$ is also be defined on this interval. We denote by $\{f(t_i^{(I)})\}_{i=1}^{N} = \{f_i^{(I)}\}_{i=1}^{N}$ the set of output values of the function $f(t)$ at the nodes of the grid $\Delta_N^{(I)}$. Let us have a look at the trigonometric polynomial

$$T_n^{(I)}(t) = \frac{a_0}{2} + \sum_{k=1}^{n} a_k^{(I)} \cos kt + b_k^{(I)} \sin kt, \qquad (1)$$

that interpolates the function $f(t)$ on the grid $\Delta_N^{(I)}$. Then the coefficients of this polynomial are determined by the formulas

$$a_k^{(I)} = \frac{2}{N} \sum_{j=1}^{N} f_j^{(I)} \cos k t_j^{(I)}, \qquad b_k^{(I)} = \frac{2}{N} \sum_{j=1}^{N} f_j^{(I)} \sin k t_j^{(I)}, \qquad (2)$$
$$k = 0, 1, ..., n; \qquad k = 1, 2, ..., n.$$

In [12], the concept of stitching grids and interpolation grids was introduced. Recall that a stitching grid is a grid in the nodes of which polynomial splines, analogs of trigonometric splines, are stitched together; we will denote the stitching grid $I1$, ($I1 = 0,1$).

An interpolation grid is a grid in the nodes of which interpolation of given values is performed by trigonometric splines or their polynomial analogs; we will denote the interpolation grid as $I2$, ($I2 = 0,1$).

Additionally, we assume that the stitching grid $I1$ and the interpolation grid $I2$ coincide with the grids $\Delta_N^{(I)}$.

Note that here and in the future, we will refer to stitching and interpolation grids as identical if $I1 = I2$ and as different if $I1 \neq I2$.

Trigonometric splines with Riemann convergence factors of degree $r$, ($r = 0,1,...$), with a stitching grid $I1$ that interpolates a given function $f$ at the nodes of the interpolation grid $I2$, were denoted as $St(f, I1, I2, r, N, t)$. However, for simplicity in this work, we will use a shortened notation for such trigonometric splines, namely $St(I1, I2, r, t)$

To illustrate the concepts of stitching and interpolation grids, let's consider an example. Let, $N = 9$ and $\{f_i\}_{i=1}^{9} = \{2,1,3,2,4,1,3,1,3\}$. We present the graphs of trigonometric interpolation splines $St(I1, I2, 1, t)$, where the stitching and interpolation grids are most prominently highlighted, along with the trigonometric interpolation polynomials $T_n^{(I2)}(t)$. Note that here and in the future, on the graphs presented below, the vertical lines coincide with the nodes of the grid $\Delta_9^{(0)}$; it is clearly that the nodes of the grid $\Delta_9^{(1)}$ are located between the nodes of the grid $\Delta_9^{(0)}$.

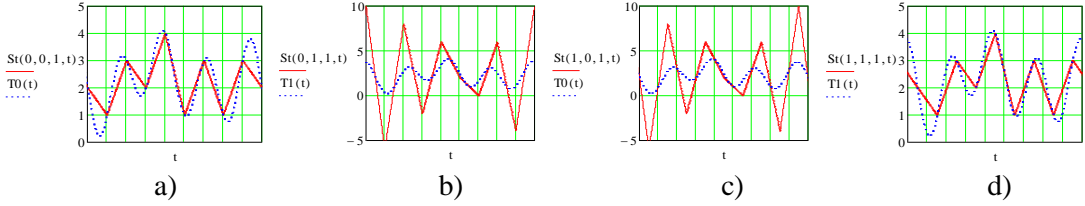

|  a) | b) | c) | d) |

Fig.1. Trigonometric splines $St(I1,I2,1,t)$:

a) $St(0,0,1,t)$ - stitching grid $I1=0$, interpolation grid $I2=0$;
b) $St(0,1,1,t)$ - stitching grid $I1=0$, interpolation grid $I2=1$;
c) $St(1,0,1,t)$ - stitching grid $I1=1$, interpolation grid $I2=0$;
d) $St(1,1,1,t)$ - stitching grid $I1=1$, interpolation grid $I2=1$.

Previously, it was shown by us that trigonometric splines $St(I1,I2,1,t)$ can be represented by the formula

$$St(I1,I2,r,t) = \frac{a_0}{2} + \qquad (3)$$

$$+\sum_{k=1}^{n} H_k^{-1}(I1,I2,r)\left(C_k(I1,I2,r,t)(a_k^{(0)})^{1-I2}(a_k^{(1)})^{I2} + S(I1,I2,r,k,t)(b_k^{(0)})^{1-I2}(b_k^{(1)})^{I2}\right),$$

where

$$H_k(I1,I2,r) = \sigma_k(r) + \sum_{m=1}^{\infty}(-1)^{m(r+1+I1+I2)}\left[\sigma_{mN+k}(r)+\sigma_{mN-k}(r)\right] \text{ - interpolation multiplier;}$$

$$C_k(I1,I2,r,t) = \sigma_k(r)\cos(kt) +$$

$$+\sum_{m=1}^{\infty}(-1)^{m(r+1+I1)}\left[\sigma_{mN+k}(r)\cos((mn+k)t)+\sigma_{mN-k}(r)\cos((mn-k)t)\right],$$

$$S_k(I1,I2,r,t) = \sigma_k(r)\sin(kt) +$$

$$+\sum_{m=1}^{\infty}(-1)^{m(r+1+I1)}\left[\sigma_{mN+k}(r)\sin((mn+k)t)-\sigma_{mN-k}(r)\sin((mn-k)t)\right],$$

$$\sigma_k(r) = \left(\frac{\sin\frac{\pi k}{N}}{\frac{\pi k}{N}}\right)^{1+r} \text{ - Riemann convergence multipliers [16].}$$

Expression (1) is a trigonometric Fourier series, the coefficients of which (denoted by $A_{mN\pm k}$ and $B_{mN\pm k}$) up to sign are products of the corresponding values of the interpolation multiplier, Riemann convergence multipliers, and coefficients of the interpolation trigonometric polynomial $T_n^{(I)}(t)$, namely

$$A_{mN\pm k}^{(I)} = H_k^{-1}(I1,I2,r)\sigma_{mN\pm k}(r)a_k^{(I)}, \quad B_{mN\pm k}^{(I)} = H_k^{-1}(I1,I2,r)\sigma_{mN\pm k}(r)b_k^{(I)}. \qquad (4)$$

It is known [17], [18] that the sum of such a series can be represented as the convolution of the sums of two series with coefficients that are the product of the coefficients of the original series. However, these two series can be obtained in several ways; below, we will consider two of them.

**1st method**. Considering that

$$\sigma_{mN\pm k}(r) = \sigma_{mN\pm k}(r-k-1)\sigma_{mN\pm k}(k), \qquad (5)$$

let us break down (2) into multipliers as follows

$$A_{mN\pm k}^{(I)} = \left[H_k^{-1}(I1,I2,r)\sigma_{mN\pm k}(0)a_k^{(I)}\right]\left[\sigma_{mN\pm k}(r-1)\right],$$

$$B_{mN\pm k}^{(I)} = \left[H_k^{-1}(I1,I2,r)\sigma_{mN\pm k}(0)b_k^{(I)}\right]\left[\sigma_{mN\pm k}(r-1)\right]. \qquad (6)$$

In the following, we need to consider the cases of even and odd values of the parameter $r$, ($r=0,1,...$).

Let us introduce functions $KR0(I1, I2, 2j, t)$ for even values of $r$ and $KR1(I1, I2, 2j-1, t)$ for odd values of $r$, which are represented by series whose coefficients are the first multipliers in formula (6):

$$KR0(I1, I2, 2j, t) = \frac{a_0}{2} + \qquad (7)$$

$$+ \sum_{k=1}^{n} H_k^{-1}(I1, I2, 2j)\left(C0_k(I1, I2, t)(a_k^{(0)})^{1-I2}(a_k^{(1)})^{I2} + S0(I1, I2, t)(b_k^{(0)})^{1-I2}(b_k^{(1)})^{I2}\right),$$

where

$$C0_k(I1, I2, t) = \sigma_k(0)\cos(kt) +$$

$$+ \sum_{m=1}^{\infty}(-1)^{m(1+I1)}\left[\sigma_{mN+k}(0)\cos((mn+k)t) + \sigma_{mN-k}(0)\cos((mn-k)t)\right],$$

$$S0_k(I1, I2, t) = \sigma_k(0)\sin(kt) +$$

$$+ \sum_{m=1}^{\infty}(-1)^{m(1+I1)}\left[\sigma_{mN+k}(0)\sin((mn+k)t) - \sigma_{mN-k}(0)\sin((mn-k)t)\right];$$

$$(j = 1, 2, \ldots).$$

$$KR1(I1, I2, 2j-1, t) = \frac{a_0}{2} + \qquad (8)$$

$$+ \sum_{k=1}^{n} H_k^{-1}(I1, I2, 2j-1)\left(C0_k(I1, I2, t)(a_k^{(0)})^{1-I2}(a_k^{(1)})^{I2} + S0(I1, I2, t)(b_k^{(0)})^{1-I2}(b_k^{(1)})^{I2}\right),$$

where

$$C1_k(I1, I2, t) = \sigma_k(0)\cos(kt) +$$

$$+ \sum_{m=1}^{\infty}(-1)^{m(I1)}\left[\sigma_{mN+k}(0)\cos((mn+k)t) + \sigma_{mN-k}(0)\cos((mn-k)t)\right],$$

$$S1_k(I1, I2, t) = \sigma_k(0)\sin(kt) +$$

$$+ \sum_{m=1}^{\infty}(-1)^{m(I1)}\left[\sigma_{mN+k}(0)\sin((mn+k)t) - \sigma_{mN-k}(0)\sin((mn-k)t)\right],$$

$$(j = 1, 2, \ldots).$$

Let's also consider the function $BR(r, t)$, which is represented by a series whose coefficients are the second multiplier in formula (6):

$$BR(r, t) = \frac{1}{\pi}\left[\frac{1}{2} + \sum_{k=1}^{n} C_k(r, t)\right], \qquad (9)$$

where

$$C_k(r, t) = \sigma_k(r)\cos(kt) + \sum_{m=1}^{\infty}\left[\sigma_{mN+k}(r)\cos((mn+k)t) + \sigma_{mN-k}(r)\cos((mn-k)t)\right].$$

It is clear that the functions $KR0(I1, I2, r, t)$ and $KR1(I1, I2, r, t)$ depend on the combination of the stitching and interpolation grids. They also contain the coefficients of the interpolation trigonometric polynomial, and thus are carriers of information about the value of the interpolated function $f(t)$ at the nodes of the interpolation grids; it is natural to call these functions kernels of interpolation trigonometric splines with Riemann convergence multipliers of the first kind.

The functions $BR(r, t)$ are the carriers of information about the smoothness of trigonometric splines and do not depend on combinations of stitching and interpolation grids; these functions on the interval $[-\pi, \pi]$ coincide with polynomial $B$-splines of degree $r$, ($r = 0, 1, \ldots$) [19]. This gives grounds to call these functions Riemann trigonometric $B$-splines of the first kind.

It is easy to see that the Riemann trigonometric $B$-splines of the first kind are normalized, i.e., they satisfy the conditions

$$\int_{-\pi}^{\pi} BR(r,t)dt = 1, \qquad (r = 0,1,...).$$

Let's consider an example that illustrates the theoretical concepts introduced. As before, let us assume $N = 9$, $\{f_i\}_{i=1}^{9} = \{2,1,3,2,4,1,3,1,3\}$. Here are the graphs of trigonometric Riemann BR-splines of the first kind for $r = 0,1,2,3$. We also give graphs of trigonometric interpolation splines $St(I1, I2, r, t)$ together with graphs of trigonometric interpolation polynomials $T_n^{(I2)}(t)$, and graphs of kernels of the first kind of these splines for even and odd values of the parameter $r$.

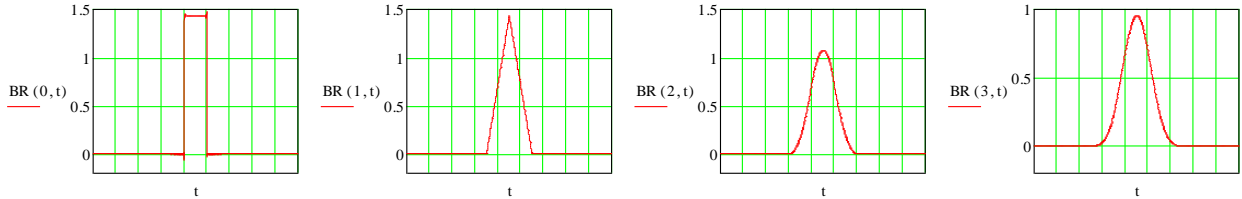

Fig.2. Trigonometric $BR(r,t)$ - splines of the first kind, ($r = 0,1,2,3$).

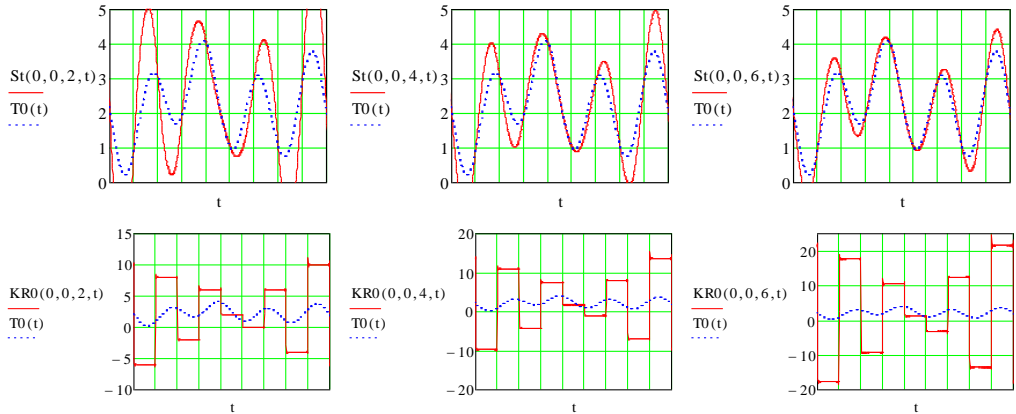

Fig.3. Trigonometric interpolation splines $St(0,0,2j,t)$
and their kernels $KR0(0,0,2j,t)$ of the first kind, ($j = 1,2,3$).

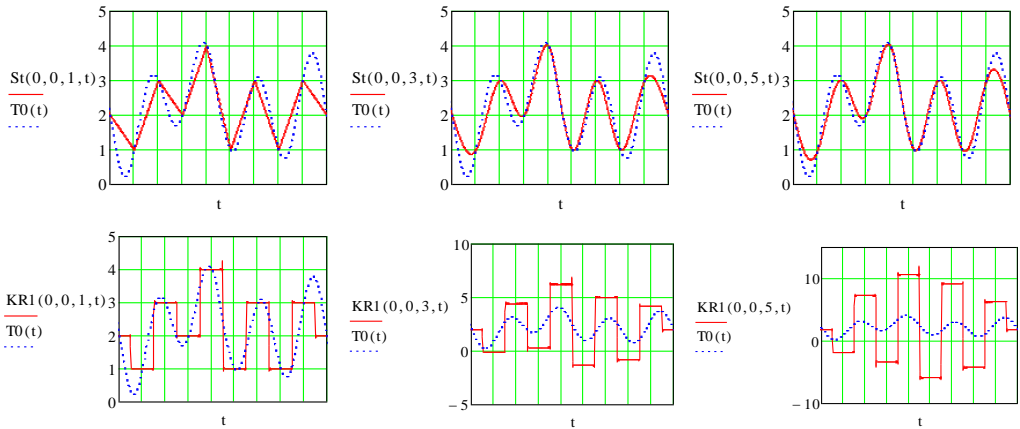

Fig.4. Trigonometric interpolation splines $St(0,0,2j-1,t)$
and their kernels $KR1(0,0,2j-1,t)$ of the first kind, ($j = 1,2,3$).

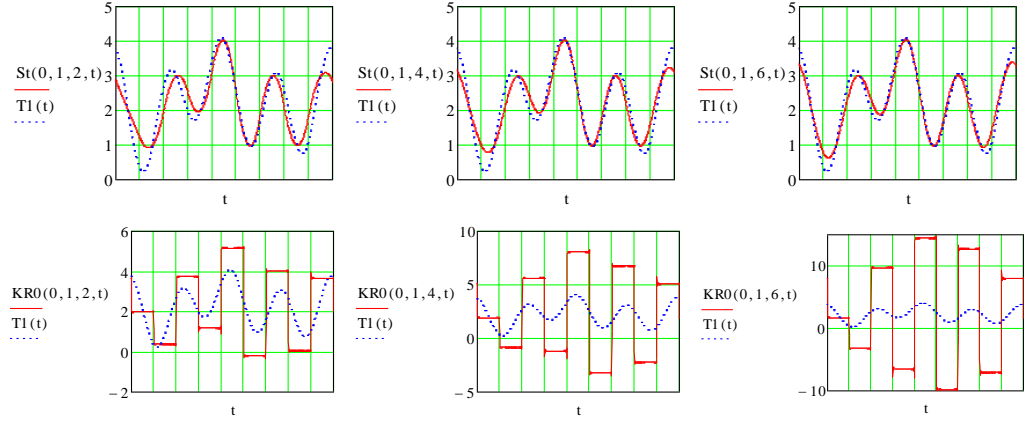

Fig.5. Trigonometric interpolation splines $St(0,1,2j,t)$

and their kernels $KR0(0,1,2j,t)$ of the first kind, ( $j = 1,2,3$ ).

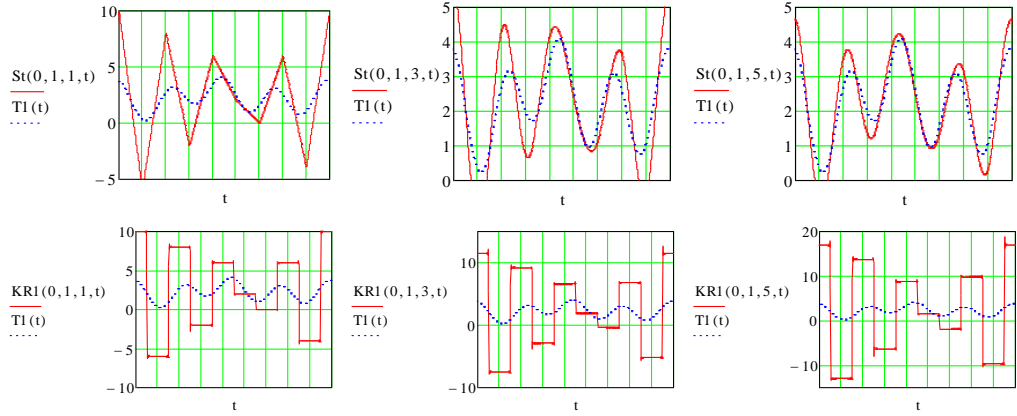

Fig.6. Trigonometric interpolation splines $St(0,1,2j-1,t)$

and their kernels $KR1(0,1,2j-1,t)$ of the first kind, ( $j = 1,2,3$ ).

We will not present the graphs of the trigonometric interpolation splines $St(1,0,r,t)$, $St(1,1,r,t)$ and their kernels of the first kind, since they do not contain any fundamentally new information.

Given the trigonometric $BR$-splines of the first kind $BR(r,t)$ and the kernels of the first kind $KR0(I1,I2,2j,t)$ and $KR1(I1,I2,2j-1,t)$, the trigonometric interpolation splines can be represented as follows:

$$St(I1,I2,2j,t) =$$
$$= \int_0^{2\pi} KR0(I1,I2,2j,t-v)BR(2j-1,v)dv = \int_0^{2\pi} KR0(I1,I2,2j,v)BR(2j-1,t-v)dv . \quad (10)$$

( $j = 1,2,...$ ).

$$St(I1,I2,2j-1,t) =$$
$$= \int_0^{2\pi} KR1(I1,I2,2j-1,t-v)BR(2(j-1),v)dv = \int_0^{2\pi} KR1(I1,I2,2j-1,v)BR(2(j-1),t-v)dv . \quad (11)$$

( $j = 1,2,...$ ).

Concluding the analysis of the first method of constructing kernels of the first kind of trigonometric interpolation $BR$-splines, we note that in this method the $BR$-splines remain constant for each fixed value of the parameter $r$, ( $r = 0,1,...$ ); however, the kernels of these splines significantly depend on the combination of stitching and interpolation grids.

**2nd method**. Let us now divide the coefficients of (4) into multipliers in a different way, namely:

$$A^{(I)}_{mN\pm k} = \left[\sigma_{mN\pm k}(0)a^{(I)}_k\right]\left[H^{-1}_k(I1,I2,r)\sigma_{mN\pm k}(r-1)\right],$$
$$B^{(I)}_{mN\pm k} = \left[\sigma_{mN\pm k}(0)b^{(I)}_k\right]\left[H^{-1}_k(I1,I2,r)\sigma_{mN\pm k}(r-1)\right]. \quad (12)$$

As before, we need to consider the cases of even and odd values of the parameter $r$, ($r = 0,1,...$). Let us introduce the functions $KR0^*(I1,I2,r,t)$ for even values of $r$ and $KR1^*(I1,I2,r,t)$ for odd values of $r$, which are represented by series whose coefficients are the first multipliers in formulas (12):

$$KR0^*(I1,I2,r,t) = \frac{a_0}{2} + \sum_{k=1}^{n}\left(C0_k(I1,I2,t)(a^{(0)}_k)^{1-I2}(a^{(1)}_k)^{I2} + S0(I1,I2,t)(b^{(0)}_k)^{1-I2}(b^{(1)}_k)^{I2}\right), \quad (13)$$

where the functions $C0_k(I1,I2,t)$ i $S0_k(I1,I2,t)$ ($k=1,2,...,n$) are the same as in formula (7).

$$KR1^*(I1,I2,r,t) = \frac{a_0}{2} + \sum_{k=1}^{n}\left(C1_k(I1,I2,t)(a^{(0)}_k)^{1-I2}(a^{(1)}_k)^{I2} + S1(I1,I2,t)(b^{(0)}_k)^{1-I2}(b^{(1)}_k)^{I2}\right), \quad (14)$$

where the functions $C1_k(I1,I2,t)$ i $S1_k(I1,I2,t)$ ($k=1,2,...,n$) are the same as in formula (8).

Let's also consider the function $BR^*(r,t)$, which is represented by a series whose coefficients are the second multiplier in formulas (12):

$$BR^*(r,t) = \frac{1}{\pi}\left[\frac{1}{2} + \sum_{k=1}^{n}\frac{C_k(r,t)}{H_k(I1,I2,1+r)}\right], \quad (15)$$

where the multiplier $H_k(I1,I2,1+r)$ is the same as in formula (3), and the function $C_k(r,t)$ is the same as in formula (9).

It is clear that now the functions $KR0^*(I1,I2,r,t)$ and $KR1^*(I1,I2,r,t)$ do not depend on the combination of stitching and interpolation grids. They also contain the coefficients of the interpolation trigonometric polynomial, and thus are carriers of information about the value of the interpolated function $f(t)$ at the nodes of the interpolation grids; as before, it is natural to call these functions kernels of the second kind of interpolation trigonometric splines with Riemann convergence multipliers.

The functions $BR^*(r,t)$ are still the carriers of information about the smoothness of trigonometric splines. However, now they depend on combinations of stitching and interpolation grids. We will refer to these functions as Riemann trigonometric $B$-splines of the second kind.

It is easy to see that Riemann trigonometric $B$-splines of the second kind are normalized, i.e., they satisfy the conditions

$$\int_{-\pi}^{\pi} BR^*(r,t)dt = 1, \qquad (r=0,1,...).$$

As before, we present graphs of trigonometric Riemann BR-splines of the second kind. We also present graphs of kernels of the second kind of these splines for even and odd values of the parameter $r$.

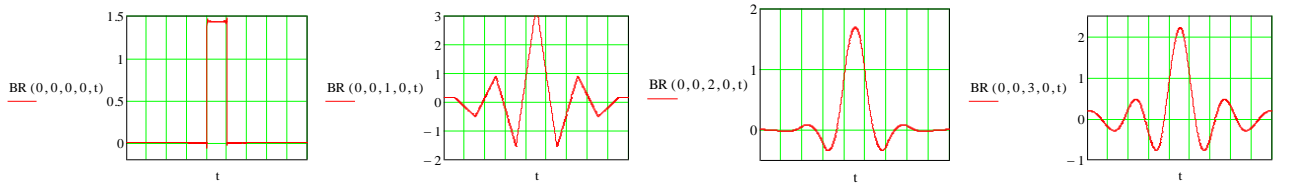

Fig.7. Riemann Trigonometric $BR^*(0,0,r,t)$ - splines of the second kind on identical grids $I1=0$, $I2=0$.

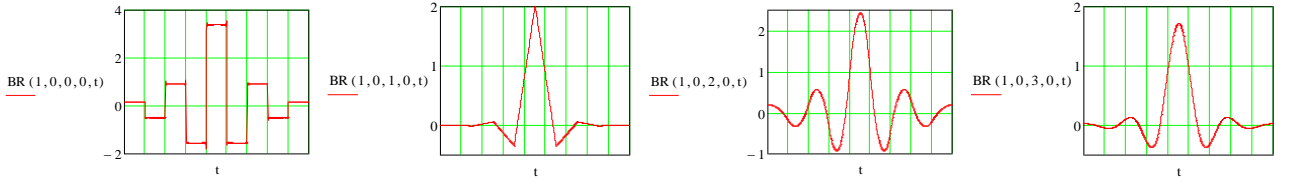

Fig.8. Riemann Trigonometric $BR^*(1,0,r,t)$ - splines of the second kind on different grids $I1 = 0, I2 = 0$.

We will not provide graphs of splines of the second kind $BR^*(0,1,r,t)$ and $BR^*(1,1,r,t)$, as the equations $BR0^*(0,0,r,t) \equiv BR0^*(1,1,r,t)$ and $BR1^*(1,0,r,t) \equiv BR1^*(0,1,r,t)$ hold.

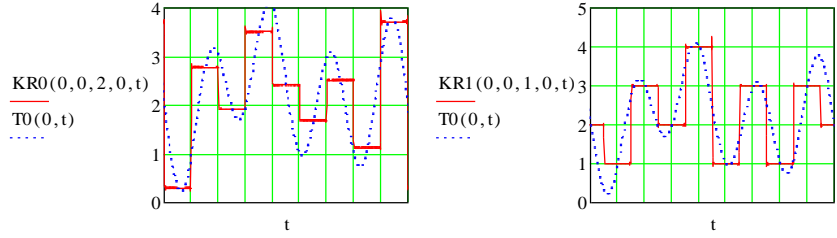

Fig.9. Kernels of the second kind $KR0^*(0,0,2j,t)$ and $KR1^*(0,0,2j-1,t)$ on identical grids.

We will not give graphs of other types of kernels of the second kind, since they satisfy the equations:

$$KR0^*(0,0,2j,t) = KR0^*(1,1,2j,t) = KR1^*(1,0,2j-1,t) = KR1^*(0,1,2j-1,t+h/2); \quad (16)$$

$$KR1^*(0,0,2j-1,t) = KR1^*(1,1,2j-1,t+h/2) = KR0^*(1,0,2j,t) = KR0^*(0,1,2j,t). \quad (16a)$$

Note that equations (16) also hold for other values of the parameter included in them; however, we will not focus on this issue here.

Given kernels of the second kind $KR0^*(I1,I2,r,t)$, $KR1^*(I1,I2,r,t)$ and $BR$-splines of the second kind $BR^*(I1,I2,r,t)$, trigonometric interpolation splines $St(I1,I2,r,t)$ can be represented as such convolutions:

$$St(I1,I2,2j,t) = \int_0^{2\pi} KR0^*(I1,I2,2j,t-v)BR^*(I1,I2,2j-1,v)dv; \quad (17)$$

$$St(I1,I2,2j-1,t) = \int_0^{2\pi} KR1^*(I1,I2,2j-1,t-v)BR^*(I1,I2,2(j-1),v)dv. \quad (17a)$$

Finally, let us formulate what we believe to be an important question. It is well known [3], [4] that polynomial B-splines are a basis in the space of polynomial splines; this property is widely used (see, e.g., [7], [8]). Trigonometric B-splines $BR(r,t)$ also form a basis in the space of trigonometric splines; this follows from the identity of periodic polynomial and trigonometric interpolation splines. It can be assumed that trigonometric splines $BR^*(r,t)$ also form a basis in the space of both trigonometric and periodic polynomial splines; however, this fact requires proof.

## Conclusions.

1. The method for constructing trigonometric B-splines and kernels of interpolation trigonometric splines with Riemann convergence multipliers has been generalized for various types of combinations of stitching and interpolation grids of trigonometric splines.

2. The kernels $KR0(I1,I2,r,t)$, $KR1(I1,I2,r,t)$, $KR0^*(I1,I2,r,t)$ and $KR1^*(I1,I2,r,t)$ of trigonometric interpolation splines are carriers of information about the interpolated function $f(t), t \in [0, 2\pi)$ (in particular, about the grid $\Delta_N$ and the value of this function on this grid).

3. The trigonometric B-splines $BR(r,t)$ and $BR^*(r,t)$ are the carriers of information about the differential properties of the spline $St(f,\sigma,r,N,t)$. Since the trigonometric B-splines $BR0(r,t)$ coincide with polynomial B-splines of the same degree on the interval $[-\pi,\pi]$ with a given grid $\Delta_N$ on it, it is natural to call the splines $BR(r,t)$ and $BR^*(r,t)$ trigonometric B-splines.

4. The kernels of the first kind $KR0(I1,I2,r,t)$ and $KR1(I1,I2,r,t)$ depend on the order of the spline $r$, ($r=0,1,...$); the kernels of the second kind $KR0^*(I1,I2,r,t)$ and $KR1^*(I1,I2,r,t)$ do not depend on the order of the spline.

5. The interpolating trigonometric splines $St(I1,I2,r,t)$ are the result of convolution of the trigonometric B-splines $BR(r,t)$ and $BR^*(r,t)$ with kernels $KR0(I1,I2,r,t)$, $KR1(I1,I2,r,t)$, $KR0^*(I1,I2,r,t)$, $KR1^*(I1,I2,r,t)$. Since it was shown earlier that the interpolation trigonometric splines $St(0,0,2j-1,t)$ coincide with periodic polynomial simple interpolation splines of the same degree ($j=1,2,...$), this result can be extended to periodic polynomial simple interpolation splines of odd degree.

6. The theoretical principles are illustrated with examples.

7. As already mentioned, the method of constructing kernels of interpolation trigonometric splines is not unique. In this paper, we considered methods for constructing kernels of order 0, which were piecewise-constant functions. However, by similar reasoning, it is possible to construct kernels of orders $r-k$ ($0 \leq k < r$).

8. Undoubtedly, the above theoretical positions require further research.